\documentstyle[twoside]{article}
\def\beq{\begin{equation}}
\def\eeq{\end{equation}}
\def\bea{\begin{eqnarray}}
\def\eea{\end{eqnarray}}
\def\nn{\nonumber}
\def\ba{\begin{array}}
\def\ea{\end{array}}
\def\one{1\hskip -0.8mm{\rm l}}
\begin{document}
\renewcommand{\thefootnote}{\fnsymbol{footnote}}
\setcounter{footnote}{1}

\begin{center}
{\bf THE EXPONENTIAL MAP FOR
REPRESENTATIONS \\ OF $U_{p,q}(gl(2))$}\footnote{Presented at the 4th
International Colloquium ``Quantum Groups and Integrable Systems,''
Prague, 22-24 June 1995.}

\vspace{1cm}

{\sc Joris Van der Jeugt}\footnote{Senior Research Associate of N.F.W.O.
(National Fund for Scientific Research of Belgium); \\
E-mail~: Joris.VanderJeugt@rug.ac.be}\\[2mm]
{\sl Department of Applied Mathematics and Computer Science,} \\
{\sl University of Ghent, Krijgslaan 281-S9, B-9000 Gent, Belgium
}\\[3mm]
{\sc Ramaswamy Jagannathan}\footnote{E-mail~: jagan@imsc.ernet.in}\\[2mm]
{\sl The Institute of Mathematical Sciences,}\\
{\sl Madras - 600113, India}
\end{center}

\vspace{1cm}

{\small For the quantum group $GL_{p,q}(2)$ and the corresponding
quantum algebra $U_{p,q}(gl(2))$ Fronsdal and Galindo~[1]
explicitly constructed the so-called universal $T$-matrix. In a
previous paper~[2] we showed how this universal $T$-matrix can be
used to exponentiate representations from the quantum algebra
to get representations (left comodules) for the quantum group. Here,
further properties of the universal $T$-matrix are illustrated. In
particular, it is shown how to obtain comodules of the quantum algebra
by exponentiating modules of the quantum group. Also the relation with
the universal $R$-matrix is discussed.}

\vspace{5mm}
\begin{center}
{\bf 1\ \ Introduction and notation}
\end{center}

Quantum groups may be studied in two different ways~: the
Drinfel'd-Jimbo approach from the point of view of integrable
theories, and the Woronowicz-Manin picture from the
point of view of pseudo-groups. Between these
two approaches there is a duality relation~[3,4].
Concretely, it means one can give dual bases for the ``quantum group''
on the one side (e.g.~$GL_q(n)$) and for the ``quantum algebra'' on the
other side (e.g.~$U_q(gl(n))$), as dual Hopf algebras. An explicit
example of such a dual basis was given by Fronsdal and
Galindo~[1] for the two-parameter quantum
group $GL_{p,q}(2)$ and
its dual $U_{p,q}(gl(2))$. This particular duality had been discussed
before~[5--7], but now Fronsdal and Galindo used it to
construct the universal $T$-matrix (also called the canonical
element~[8]) yielding the interpretation
as an exponential mapping from the quantum algebra to the quantum group.
This work has been extended to other cases~[9--11].
Here and
in~[2] we show that this exponential mapping can be applied to
representations. In particular, we show how left modules of the
quantum algebra give rise to left comodules of the quantum group. Using the
dual pairing this also leads to left modules of the quantum group and
the related left comodules of the quantum algebra.

For the quantum group $GL_{p,q}(2)$ the defining $T$-matrix, or the defining
representation matrix, is specified by
\beq
T = \left( \ba{cc}
a & b \\
c & d \\
\ea \right)
\label{defT}
\eeq
with the commutation relations in terms of two parameters $p$ and $q$~:
\bea
ab & = & qba, \qquad cd = qdc, \qquad ac = pca, \qquad bd = p db, \nn \\
bc & = & (p/q)cb, \qquad ad-da = (q-p^{-1})bc,
\label{comrel}
\eea
following from the $RTT$-relation $RT_1T_2=T_2T_1R$, with
$T_1=T\otimes\one$, $T_2=\one\otimes T$ and
\beq
R = Q^{\frac{1}{2}}\left(
\ba{cccc}
Q^{-1} & 0 & 0 & 0 \\
0 & \lambda^{-1} & Q^{-1}-Q & 0 \\
0 & 0 & \lambda & 0 \\
0 & 0 & 0 & Q^{-1}
\ea \right),
\quad Q = \sqrt{pq}, \qquad  \lambda = \sqrt{p/q}\,.
\label{R}
\eeq
The coproduct and counit are easily defined in terms of $T$, e.g.\
$\Delta(T)=T\dot\otimes T$, and an antipode can be defined making
$A_{p,q}(GL(2))$, the algebra of functions on $GL_{p,q}(2)$, into a Hopf
algebra.
For the definition of the antipode, ${\cal D} = ad - qbc = ad - pcb$,
the quantum determinant of $T$, satisfying
$\Delta ({\cal D})$ $=$ ${\cal D} \otimes {\cal D}$, should
be invertible.  The noncentral ${\cal D}$ satisfies
\beq
{\cal D}a = a{\cal D}, \qquad {\cal D}b = \lambda^{-2}b{\cal D}, \qquad
{\cal D}c = \lambda^2 c{\cal D}, \qquad {\cal D}d = d{\cal D} .
\label{qdcomrel}
\eeq

The algebra ${\cal U}$ dual to ${\cal A}$ $=$ $A_{p,q}(GL(2))$,
namely $U_{p,q}(gl(2))$, is generated by the elements
$\{ J_0, J_\pm\,,Z \}$ subject to the relations
\beq
[ J_0,J_\pm ] =  \pm J_\pm, \quad \quad
[ J_+,J_- ] = [ 2J_0 ]_Q,   \quad \quad
{} [ Z , \cdot\,] = 0,
\label{alg}
\eeq
where, as usual, $[x]_t=(t^x-t^{-x})/(t-t^{-1})$. We recall here the
form of the coproduct~:
\bea
\Delta ( J_\pm ) & =
& J_\pm \otimes Q^{-J_0}\lambda^{\pm Z} +
Q^{J_0}\lambda^{\mp Z} \otimes J_\pm, \nn \\
\Delta ( J_0 ) & =
& J_0 \otimes \one + \one \otimes J_0, \qquad
\Delta (Z) = Z \otimes \one + \one \otimes Z;
\label{coalg}
\eea
the remaining functions (counit, antipode) can be found
elsewhere~[6,7].
Apart from $\Delta$, there exists also
the opposite coproduct $\Delta '=\sigma\Delta$, with
$\sigma(u \otimes v)=v \otimes u$. For a quasitriangular Hopf algebra
such as ${\cal U}=U_{p,q}(gl(2))$ these two coproducts are related
by $\Delta'{\cal R} = {\cal R}\Delta $
through an
element ${\cal R}\in{\cal U}\otimes{\cal U}$ called the universal
$R$-matrix and explicitly given by~[12]
\bea
{\cal R} & = & Q^{-2(J_0 \otimes J_0)}
\lambda^{2(Z \otimes J_0 - J_0 \otimes Z)} \nn \\
&    &  \ \ \times \sum_{n=0}^{\infty}\,
\frac{(1-Q^2)^n}{[n]!}\,Q^{-\frac{1}{2}
n(n-1)}(Q^{-J_0}\lambda^Z J_+
\otimes Q^{J_0}\lambda^Z J_-)^n
\label{calR}
\eea
where $[n] = [n]_Q$ and $[n]! = [n][n-1]\ldots[2][1]$.

\begin{center}
{\bf 2\ \ The universal $T$-matrix}
\end{center}

The exponential map from ${\cal U}=U_{p,q}(gl(2))$ to
$A_{p,q}(GL(2))$, {\em \`{a} la} Fronsdal and
Galindo~[1]
is defined in terms of the so-called universal $T$-matrix ${\cal T}$.
First, one has to construct a dual basis for the dual Hopf algebras
${\cal A}$ and ${\cal U}$. The fact that ${\cal A}$ and ${\cal U}$ are
dual means that product, coproduct, unit, counit and antipode in ${\cal
A}$ and ${\cal U}$ are related as follows by the dual pairing~:
\bea
&& \langle XY,x\rangle=\langle X\otimes Y, \Delta x\rangle, \qquad
\langle \Delta X,x\otimes y\rangle=\langle X, xy\rangle, \qquad
\langle 1,x \rangle = \epsilon(x),\nn\\
&& \qquad \epsilon(X)=\langle
X,1\rangle, \qquad \langle S X,x\rangle = \langle X, S x
\rangle,\qquad X,Y\in{\cal U},\ x,y\in{\cal A}.
\eea
A dual basis then consists of a basis
$\{ x^A \}$ for ${\cal A}$ and $\{ X_A\}$ for ${\cal U}$ such that the
pairing is equal to $\langle X_A, x^B \rangle = \delta_A^B$ (Kronecker
$\delta$). As a consequence, the structure constants for the product
(resp.\ coproduct) in ${\cal A}$ become the structure constants for the
coproduct (resp.\ product) in ${\cal U}$~:
\bea
&& x^Bx^C = \sum_{A}\,h^{BC}_A x^A, \quad \quad
\Delta (x^C) = \sum_{A,B}\,f^C_{AB} x^A \otimes x^B , \nn\\
&& X_AX_B = \sum_{C}\,f_{AB}^C X_C, \quad \quad
\Delta (X_A) = \sum_{B,C}\,h_A^{BC} X_B \otimes X_C\,.
\label{XX}
\eea
The universal ${T}$-matrix, ${\cal T}$, can now be defined as an
element of ${\cal A}\otimes {\cal U}$~:
\beq
{\cal T} = \sum_{A}\,x^A\otimes X_A.
\label{TauxX}
\eeq
Using the duality properties, it is easy to show that ${\cal T}$
satisfies the following properties~:
\bea
&& (\Delta\otimes\one) ({\cal T}) = \sum_{A}\,\Delta (x^A)\otimes X_A
= {\cal T}_{13} {\cal T}_{23},\nn\\
&& (\one\otimes\Delta) ({\cal T}) = \sum_{A}\,x^A\otimes \Delta(X_A)
= {\cal T}_{12} {\cal T}_{13},
\label{DelcalTau}
\eea
where, e.g.\ ${\cal T}_{13}$ stands for $\sum_A x^A\otimes 1\otimes X_A$.

What is now this dual basis and the corresponding universal $T$-matrix
for the case ${\cal U}=U_{p,q}(gl(2))$ and
${\cal A}=A_{p,q}(GL(2))$? To this end, the generators of
$U_{p,q}(gl(2))$ are redefined as follows~:
\beq
\hat{J}_+ = J_+Q^{-(J_0+\frac{1}{2})}
\lambda^{Z-\frac{1}{2}}, \qquad
\hat{J}_- = Q^{(J_0+\frac{1}{2})}
\lambda^{Z-\frac{1}{2}}J_-, \qquad
\hat{J}_0 = J_0,\qquad
\hat{Z} = Z  ,
\label{newJ}
\eeq
with, of course, the corresponding products and coproducts. For
${\cal A}$ $=$ $A_{p,q}(GL(2))$, with the extra assumption
that $a$ is invertible, the new variables $\{\alpha ,\beta ,
\gamma ,\delta \}$ defined through
\beq
a=e^\alpha,\qquad \beta = a^{-1}b, \qquad \gamma = ca^{-1}, \qquad
d=ca^{-1}b + e^{-\delta}.
\label{begadel}
\eeq
form a Lie algebra with structure constants depending upon
$\rho=\ln(Q)$ and $\theta=\ln(\lambda)$~:
\bea
[ \alpha , \beta ] & = & (\rho -\theta )\beta, \quad \quad
[ \alpha , \gamma ] = (\rho + \theta )\gamma, \nn \\
{} [ \delta , \beta ] & = & (\rho + \theta )\beta, \quad \quad
[ \delta , \gamma ] = (\rho - \theta )\gamma, \nn \\
{} [ \alpha , \delta ] & = & 0, \quad \quad
[ \beta , \gamma ] = 0 .
\label{Lie}
\eea
With these new generating elements, Fronsdal and Galindo~[1]
obtained an explicit dual basis for $A_{p,q}(GL(2))$ and
$U_{p,q}(gl(2))$ given by~:
\bea
x^A & = & \gamma^{a_1}\alpha^{a_2}\delta^{a_3}\beta^{a_4}, \nn \\
X_A & = &
\frac{Q^{\frac{1}{2}a_1(a_1-1)}\hat{J}_-^{a_1}}
{[a_1]!}
\frac{(\hat{J}_0 + \hat{Z})^{a_2}}{a_2!}
\frac{(\hat{J}_0 - \hat{Z})^{a_3}}{a_3!}
\frac{Q^{-\frac{1}{2}a_4(a_4-1)}\hat{J}_+^{a_4}}
{[a_4]!}, \nn \\
    &   &   a_1,a_2,a_3,a_4 = 0,1,2,\ldots\,.
\label{defxX}
\eea
Then, with the definition of a basic exponential function
with parameter $t^2$ as
\beq
{\cal E}_{t^2}(x) =
\sum_{n=0}^{\infty}\,
\frac{t^{-\frac{1}{2}n(n-1)}}{[n]_t!}\,x^n,
\label{calE}
\eeq
the universal $T$-matrix~(\ref{TauxX}) can be written as
\beq
{\cal T} = {\cal E}_{Q^{-2}}(\gamma\otimes \hat{J}_-)\;
{\rm e}^{\alpha\otimes(\hat{J}_0+\hat{Z}) +
\delta\otimes(\hat{J}_0-\hat{Z})}\;
{\cal E}_{Q^2}(\beta\otimes \hat{J}_+).
\label{calT}
\eeq

\begin{center}
{\bf 3\ \ Applications}
\end{center}

The first application was discussed in~[2]. There, it was shown
that (\ref{calT}) can be used to exponentiate representations of
$U_{p,q}(gl(2))$ to representations of $GL_{p,q}(2)$. For generic
parameters, the finite dimensional irreducible representations of
$U_{p,q}(gl(2))$ are labelled by two numbers, the `spin' $j$
(integer or half-integer) and the $Z$-eigenvalue $z$. The
explicit
matrix elements of the generating elements in this $(2j+1)\times (2j+1)$
representation $(j;z)$ are given by
\bea
(J_\pm)_{mk} & = & ([j \pm m][j+1 \mp m])^{1/2} \delta_{m,k\pm 1}, \nn \\
(J_0)_{mk} & = & m\delta_{mk}, \qquad
(Z)_{mk} = z\delta_{mk}, \nn \\
  &   &  \quad m,k = j,j-1,\ldots ,-(j-1),-j\,.
\label{Gamma}
\eea
Exponentiating this representation means substituting these matrices in
the formula (\ref{calT}). This gives rise to a matrix of order $2j+1$
with elements depending upon $\alpha$, $\beta$, $\gamma$ and $\delta$,
or using the relations (\ref{begadel}) eventually upon $a,b,c$ and $d$.
This matrix, denoted by $T^{(j;z)}$, was calculated explicitly
in~[2], and its elements are given as follows~:
\bea
T^{(j;z)}_{mk} & = & {\cal D}^{z-j}
Q^{-(m-k)(2j-m+k)/2}\lambda^{-(m-k)(2j-2z-m-k)} \nn\\
 &\times &   ([j+m]![j-m]![j+k]![j-k]!)^{1/2}  \nn \\
   & \times &\sum_{s}\,Q^{-s(2j-m+k-s)}\lambda^{-s(m-k+s)} \nn\\
  &\times &    \frac{a^{j+k-s}}{[j+k-s]!} \frac{b^{m-k+s}}{[m-k+s]!}
                \frac{c^s}{[s]!} \frac{d^{j-m-s}}{[j-m-s]!}\,.
\label{Taumk}
\eea
In the limit $p = q$, or $Q = q$ and $\lambda = 1$, and ${\cal D}=1$,
$GL_{p,q}(2)$ becomes $SL_q(2)$. Then, the above matrices coincide with
the representations of $SL_q(2)$ as obtained earlier in different
approaches~[13--19].

As an example, for $(j;z)=(1/2;1/2)$, the matrix $T^{(j;z)}$ simply
becomes the defining $T$-matrix given in (\ref{defT}). For $(j;z)=(1;1/2)$,
one finds
\beq
T^{(1;1/2)} = {\cal D}^{-1/2} \left(
\ba{ccc}
a^2 &
[2]^{\frac{1}{2}}Q^{-\frac{1}{2}}ab &
\lambda^{-1}b^2 \\
  &   &  \\
{} [2]^{\frac{1}{2}}Q^{-\frac{1}{2}}ac &
ad+Q^{-1}\lambda^{-1}bc &
[2]^{\frac{1}{2}}Q^{-\frac{1}{2}}\lambda^{-1}bd \\
  &   &  \\
\lambda c^2 &
[2]^{\frac{1}{2}}Q^{-\frac{1}{2}}\lambda cd &
d^2
\ea \right)\,.
\label{Tau1abcd}
\eeq

Because of the properties of the universal $T$-matrix ${\cal T}$, the
matrices $T^{(j;z)}$ satisfy the following properties~:
\beq
\Delta(T^{(j;z)}_{kl})=\sum_{i=-j}^j T^{(j;z)}_{ki}\otimes
T^{(j;z)}_{il}, \qquad {\rm i.e.}\qquad
\Delta(T^{(j;z)})=T^{(j;z)}\dot\otimes T^{(j;z)}.
\eeq
In particular, if the elements of a column of the matrix $T^{(j;z)}$
are denoted by $v_l$, then $\Delta(v_l)=\sum_i T^{(j;z)}_{li}\otimes
v_i$. This means that we are dealing with a left ${\cal A}$-comodule.

For the second application, we shall take ${\cal A}=A_q(SL(2))$
and ${\cal U}=U_q(sl(2))$ in which case $p=q$
(or, $Q=q$, $\lambda=1$) and ${\cal D}=1$.
It is well
known that the relations (and other properties such as $\Delta$,
$\epsilon$ and
$S$) for ${\cal U}$ can also be described in terms of the
standard $L$-matrices $L^\pm$, together with the standard $R$-matrix
(\ref{R})~:
\beq
RL_2^\pm L_1^\pm=L_1^\pm L_2^\pm R,\qquad
RL_2^+ L_1^- = L_1^- L_2^+ R.
\eeq
Explicitly, we have
\beq
L^+=\left( \begin{array}{cc} q^{-J_0} & 0 \\ q^{1/2}(q^{-1}-q)J_+ &
q^{J_0} \end{array}\right) ,\qquad
L^-=\left( \begin{array}{cc} q^{J_0} & q^{-1/2}(q-q^{-1})J_- \\ 0 &
q^{-J_0} \end{array}\right) .
\eeq
It is known that these matrices satisfy the defining relations for
$SL_{q^{-1}}(2)$.  Hence, changing $q$ to $q^{-1}$ in the above
expressions leads to two maps $\pi^-$ and $\pi^+$ which are Hopf
algebra homomorphisms~:
\beq
\pi^- : \left\{
\begin{array}{ccl} a & \rightarrow & q^{J_0} \\
 b & \rightarrow & 0 \\   c & \rightarrow & q^{-1/2}(q-q^{-1})J_+ \\
 d & \rightarrow & q^{-J_0} \end{array}\right. \qquad\qquad
\pi^+ : \left\{
\begin{array}{ccl} a & \rightarrow & q^{-J_0} \\
 b & \rightarrow & q^{1/2}(q^{-1}-q)J_- \\   c & \rightarrow & 0 \\
 d & \rightarrow & q^{J_0} \end{array}\right. .
\eeq
Given the $(2j+1)$-dimensional representation $\Gamma^{(j)}$ of
$U_q(sl(2))$ (see (\ref{Gamma}) with $Q=q$ and $\lambda=1$), applying
$\Gamma^{(j)} \circ \pi^\pm$ yields $(2j+1)$-dimensional
representations for $a,b,c,d$ (left modules for $SL_q(2)$), and thus
also matrix representations for $\alpha, \beta, \gamma, \delta$. These
representations for $a,b,c,d$ coincide with the ones obtained from the
universal $R$-matrix as follows. Let $R^{(1/2)\otimes(j)}$ be the
matrix obtained from the universal $R$-matrix by restriction to the
representations specified, $(R^{(1/2)\otimes(j)})^{-1}$ its
inverse, and $T=(t_{ik})$ the standard $T$-matrix (\ref{defT}). Then
\beq
(\Gamma^{(j)} \circ \pi^+) (t_{ik})_{lm}=R^{(1/2)\otimes(j)}_{il,km},
\qquad
(\Gamma^{(j)} \circ \pi^-)
(t_{ik})_{lm}=(R^{(1/2)\otimes(j)})^{-1}_{li,mk}.
\eeq
Substituting these representations in ${\cal T}$ as follows
\beq
L^{+(j)}= \left( (\Gamma^{(j)}\circ \pi^+)\otimes\one \right) ({\cal
T}),\qquad
L^{-(j)}= \left( (\Gamma^{(j)}\circ \pi^-)\otimes\one \right) ({\cal T}),
\eeq
one obtains $(2j+1)\times(2j+1)$ matrices with elements that are
polynomials of degree $2j$ in $J_\pm$ and $q^{\pm J_0}$. In
particular, one finds that $L^{\pm(1/2)} = L^\pm$.  For $j = 1,$
\beq
L^{+(1)}=
\left(\begin{array}{ccc}
q^{-2J_0} & 0 & 0\\
(1-q^2)\sqrt{[2]}q^{-J_0}J_+ & 1 & 0\\
q^{-1}(1-q^2)^2J_+^2 & q^{-1}(1-q^2)\sqrt{[2]}q^{J_0}J_+ & q^{2J_0}
\end{array}\right),
\eeq
\beq
L^{-(1)}=
\left(\begin{array}{ccc}
q^{2J_0} & q(1-q^{-2})\sqrt{[2]}q^{J_0}J_- & q(1-q^{-2})^2 J_-^2 \\
0 & 1 & (1-q^{-2})\sqrt{[2]}q^{-J_0}J_-\\
0 & 0 & q^{-2J_0}
\end{array}\right).
\eeq
The matrices $L^{\pm(j)}$ are the spin $j$ analogs of the standard
$L$-matrices. They satisfy the properties
\bea
&&R^{(j)\otimes(j)}L_2^{\pm(j)} L_1^{\pm(j)}=L_1^{\pm(j)} L_2^{\pm(j)}
R^{(j) \otimes(j)},\nn\\
&&R^{(j)\otimes(j)}L_2^{+(j)} L_1^{-(j)}=L_1^{-(j)} L_2^{+(j)}
R^{(j) \otimes(j)},\\
&&\Delta(L^{\pm(j)}) = L^{\pm(j)} \dot\otimes L^{\pm(j)},
\eea
and thus they constitute the $(2j+1)$-dimensional comodules of
$U_q(sl(2))$.

As a third and last application, we write ${\cal T}$ in a different
form~:
\beq
{\cal T}' = \sum_A X_A \otimes x^A \in {\cal U} \otimes {\cal A}.
\eeq
Then
\beq
(\one\otimes\pi^+ ){\cal T}' = {\cal R},\qquad
(\one\otimes\pi^- ){\cal T}' = {\cal R}.
\eeq
This relationship between the univeral $T$-matrix and the universal
$R$-matrix has also been considered by Fronsdal~[20].
\vskip 2mm
{\small One of the authors (JVdJ) would like to thank the organizers of
this colloquium for giving him the opportunity to attend this meeting.
This research was partly supported by the EU (contract CI1*-CT92-0101).}

\begin{center}
{\bf References}
\end{center}

\begin{enumerate}

\item
Fronsdal C.\ and Galindo A.: Lett. Math. Phys. {\it 27}
(1993) 59.

\item
Jagannathan R.\ and Van der Jeugt J.: J. Phys. A {\it 28} (1995) 2819.

\item
Reshetikhin N.~Yu., Takhtajan L.~A. and Faddeev L.~D.: Leningrad
Math. J. {\it 1} (1990) 193.

\item
Majid S.: Int. J. Mod. Phys. A {\it 5} (1990) 1.

\item
Schirrmacher A., Wess J. and Zumino B.: Z. Phys. C {\it 49} (1991) 317.

\item
Dobrev V.~K.: Introduction to Quantum Groups (Preprint,
G\"{o}ttingen Univ., 1991) (to appear in
Proc. 22nd Ann. Iranian Math. Conf., Mashhad, March 1991).

\item
Dobrev V.~K.: J. Math. Phys. {\it 33} (1992) 3419.

\item
Woronowicz S.L.: Comm. Math. Phys. {\it 111} (1987) 613.

\item
Bonechi F., Celeghini E., Giachetti R., Perena C.~M., Sorace E. and
Tarlini M.: J. Phys. A {\it 27} (1994) 1307.

\item
Morozov A. and Vinet L.: Free-field representation of group element for
simple quantum groups (Preprint CRM-2202, Univ. Montreal, 1994;
hep-th/9409093).

\item
Chakrabarti R.\ and Jagannathan R.:
Lett. Math. Phys. (in press, 1995).

\item
Chakrabarti R.\ and Jagannathan R.: J. Phys. A {\it 27} (1994) 2023.

\item
Woronowicz S.~L.: Publ. RIMS, Kyoto Univ. {\it 23} (1987) 117.

\item
Vaksman L.~L. and Soibel'man Ya.~S.: Func. Anal. Appl. {\it 22} (1988) 170.

\item
Koornwinder T.~H.: Indag. Math. {\it 51} (1989) 97.

\item
Groza V.~A., Kachurik I.~I. and Klimyk A.~U.: J. Math. Phys. {\it 31}
(1990) 2769.

\item
Masuda T., Mimachi K., Nakagami Y., Noumi M. and Ueno K.: J. Func.
Anal. {\it 99} (1991) 357.

\item
Nomura M.: J. Phys. Soc. Japan {\it 59} (1990) 4260.

\item
Koelink H.~T. and Koornwinder T.~H.:  Nederl. Akad. Wetensch. Proc.
Ser. A {\it 92} (1989) 443.

\item
Fronsdal C.: Universal $T$-matrix for twisted quantum $gl(N)$, preprint
UCLA\-93/TEP/3 (1993) (to appear in Proc. Nato Conf. on Quantum
Groups, San Antonio, Texas, 1993).

\end{enumerate}
\end{document}